\input amstex
\magnification 1200
\TagsOnRight
\NoBlackBoxes
\def\qed{\ifhmode\unskip\nobreak\fi\ifmmode\ifinner
\else\hskip5pt\fi\fi\hbox{\hskip5pt\vrule width4pt height6pt
depth1.5pt\hskip1pt}}
\baselineskip 20 pt
\parskip 6 pt

\centerline {\bf TRANSMISSION EIGENVALUES FOR THE SELFADJOINT}
\centerline {\bf SCHR\"ODINGER OPERATOR ON THE HALF LINE}

\vskip 10 pt
\centerline {Tuncay Aktosun}
\vskip -8 pt
\centerline{Department of Mathematics}
\vskip -8 pt
\centerline {University of Texas at Arlington}
\vskip -8 pt
\centerline{Arlington, TX 76019-0408, USA}
\vskip -8 pt
\centerline {aktosun\@uta.edu}

\centerline {Vassilis G. Papanicolaou}
\vskip -8 pt
\centerline {Department of Mathematics}
\vskip -8 pt
\centerline {National Technical University of Athens}
\vskip -8 pt
\centerline {Zografou Campus}
\vskip -8 pt
\centerline {157 80, Athens, Greece}
\vskip -8 pt
\centerline {papanico\@math.ntua.gr}

\vskip 10 pt

\noindent {\bf Abstract}:
The transmission eigenvalues corresponding to the
half-line Schr\"odinger equation with the general selfadjoint
boundary condition is analyzed when the potential is real
valued, integrable, and compactly supported. It is shown that a transmission eigenvalue
corresponds to the energy at which the scattering from the perturbed
system agrees with the scattering from the unperturbed system.
A corresponding inverse problem for the recovery of the
potential from a set containing the boundary condition and the
transmission eigenvalues is analyzed, and a unique reconstruction of the
potential is
given provided one additional constant is contained in the data set.
The results are illustrated with various explicit examples.

\vskip 3 pt
\noindent {\bf PACS (2010):}  02.30.Zz, 03.65.Nk\hfil
\vskip -8 pt
\par \noindent {\bf Mathematics Subject Classification (2010):}
34A55, 34B07, 34B24
\vskip -8 pt
\par\noindent {\bf Keywords:} Transmission eigenvalues,
inverse problem, Schr\"odinger equation on the half line,
selfadjoint boundary condition
\vskip -8 pt
\par\noindent {\bf Short title:}
Transmission eigenvalues for selfadjoint problems

\newpage

\noindent {\bf 1. INTRODUCTION}
\vskip 3 pt

We consider the so-called transmission eigenvalue problem for the
half-line Schr\"odinger operator with the general selfadjoint boundary condition at the origin. We analyze the corresponding direct and inverse problems when the potential $V$ in the Schr\"odinger equation is real valued, vanishes when $x>b$ for some positive $b,$ and integrable on the interval $(0,b).$ We say that $V$ belongs to class $\Cal A$ if it satisfies the aforementioned three conditions.
The real-valuedness and integrability are standard assumptions [7,15,16,21] on the potential
of the Schr\"odinger equation, and the compact-support property naturally arises
in the analysis of transmission eigenvalues [8,10,11]. Thus, it is reasonable to
restrict our analysis to potentials in class $\Cal A.$
Our direct problem consists of the determination of the transmission eigenvalues when the potential and the boundary condition are known.
Our inverse problem consists of the recovery of the potential from an appropriate data set
containing the transmission eigenvalues.

There are two primary reasons for us to use a general selfadjoint boundary condition at the origin rather than the Dirichlet boundary condition [7,21]. First, the use of a general selfadjoint boundary condition truly clarifies the meaning and physical interpretation of the transmission eigenvalues. Second,
there are important physical problems where selfadjoint boundary conditions
other than a Dirichlet boundary condition naturally arise. Hence, our work contributes
to the analysis of direct and inverse problems associated with transmission eigenvalues,
perhaps by being the first study to consider a general selfadjoint boundary condition instead of the
mere Dirichlet boundary condition.

Due to the presence of a boundary parameter in the
non-Dirichlet case, the analysis of the Schr\"odinger equation with non-Dirichlet boundary conditions
is naturally more elaborate than the analysis under a Dirichlet boundary condition.
There are both similarities and differences between the Dirichlet and non-Dirichlet cases. We refer the reader to [5,15,16] and the references therein for the contrast between
those cases in the analysis of (1.1). In our study of transmission eigenvalues, we mainly concentrate on the non-Dirichlet case, but we also provide in Section 7 a summary of the corresponding results
in the Dirichlet case in order to have a comparison with the non-Dirichlet case.

Thus, we consider the Schr\"odinger equation on the half line
$$-\psi'' +V(x)\,\psi=k^2\psi,\qquad x\in{\bold R^+},\tag 1.1$$
where ${\bold R^+}:=(0,+\infty),$
the prime denotes the $x$-derivative, and
the potential $V$ belongs to class $\Cal A$ and thus vanishes for $x>b.$
The most general selfadjoint boundary condition at $x=0$ associated
with (1.1) is given by [5,15,16]
$$(\sin\theta)\, \psi'(0)+(\cos \theta)\, \psi(0)=0,\tag 1.2$$
where the boundary parameter $\theta$ can take any value in the interval
$(0,\pi].$ The case $\theta=\pi$ corresponds to the Dirichlet boundary condition
and is equivalent to
$$\psi(0)=0.\tag 1.3$$
In the non-Dirichlet case, i.e. when $\theta\in (0,\pi),$ we can write (1.2) as
$$\psi'(0)+(\cot\theta)\,\psi(0)=0,\qquad 0<\theta<\pi.\tag 1.4$$
Note that the mapping $\theta\mapsto \cot\theta$ is one-to-one and onto from the interval
$(0,\pi)$ to the entire real axis ${\bold R},$ and hence (1.4) can be used for many
physical problems with
an appropriate choice of $\theta$ in the interval $(0,\pi).$

If (1.1) comes from the three-dimensional Schr\"odinger equation with a spherically
symmetric potential, then it is natural to impose (1.3) so that the corresponding
solution to the three-dimensional Schr\"odinger equation remains finite at $x=0.$
Because (1.3) is used as the implicit
boundary condition in many physical problems, some physicists may not even be aware of the
mathematical necessity of imposing a boundary condition at $x=0$ for (1.1). However,
the so-called bound-state energies corresponding to the discrete
eigenvalues of (1.1) are directly affected by the
choice of the boundary parameter $\theta$ appearing in (1.2).
We refer
the reader to (1.4)-(1.6) of [3] for the elaboration on the natural occurrence
of (1.3). On the other hand, there are important physical problems where
(1.4) rather than (1.3) is appropriate to use. For example, in the inverse problem
of the recovery of the shape of the human vocal tract from sound pressure measurements
at the mouth, (1.1) and (1.4) arise in a natural manner [1,2] with
$$V(x)=\displaystyle\frac{r''(x)}{r(x)},\quad \cot\theta=-\displaystyle\frac{r'(0)}{r(0)},$$
where $r(x)$ corresponds to the cross sectional radius of the
vocal tract as a function of the distance from the glottis,
and $r'(x)$ corresponds to the slope (bending) of
that radius function, with the understanding that
 $x=0$ indicates the location of the glottis. The boundary condition (1.4) also appears in various
other vibrating
systems [13].

The transmission eigenvalues [3,4,8-11,16-19]
for the Schr\"odinger equation with the
Dirichlet boundary condition (1.3) correspond to those $\lambda$-values yielding
nontrivial solutions $\psi$ and $\psi_0$ for the
system
$$\cases -\psi'' +V(x)\,\psi=\lambda\psi,\qquad 0<x<b,\\
\noalign{\medskip}
-\psi_0''=\lambda\psi_0,\qquad 0<x<b,\\
\noalign{\medskip}
\psi(0)=\psi_0(0)=0,\\
\noalign{\medskip}
\psi_0(b)=\psi(b),\quad \psi_0'(b)=\psi'(b).
\endcases$$
On the other hand, the transmission eigenvalues for the Schr\"odinger equation with the
non-Dirichlet boundary condition (1.4) correspond to those $\lambda$-values yielding
nontrivial solutions $\psi$ and $\psi_0$ for the
system
$$\cases -\psi'' +V(x)\,\psi=\lambda\psi,\qquad 0<x<b,\\
\noalign{\medskip}
-\psi_0''=\lambda\psi_0,\qquad 0<x<b,\\
\noalign{\medskip}
\psi'(0)+(\cot \theta)\, \psi(0)=0,\\
\noalign{\medskip}
\psi_0'(0)+(\cot \theta)\, \psi_0(0)=0,\\
\noalign{\medskip}
\psi_0(b)=\psi(b),\quad \psi_0'(b)=\psi'(b),
\endcases\tag 1.5$$
which is obtained by replacing the Dirichlet boundary condition at $x=0$
with the general selfadjoint boundary condition given in (1.4).

Our paper is organized as follows. We first analyze the direct problem
for (1.1) with the boundary condition (1.4) corresponding to the
non-Dirichlet case. Our direct problem under study consists of the determination of the
corresponding transmission eigenvalues when the potential $V$ in class $\Cal A$
and the boundary parameter $\cot\theta$ are given. For this purpose,
in Section 2, we introduce the corresponding
Jost solution $f(k,x),$ the regular solution
$\varphi(k,x),$ the Jost function $F(k),$ and the scattering matrix $S(k),$
and we present their properties relevant to our study. In Section~2, we also
introduce the quantities corresponding to (1.1) with $V(x)\equiv 0$ and (1.4),
namely the Jost solution $f_0(k,x),$ the regular solution $\varphi_0(k,x),$
the Jost function $F_0(k),$ and the scattering matrix $S_0(k),$ which are
denoted by using the subscript zero. In the same
section we indicate that a potential $V$ in class $\Cal A$
is uniquely determined by the corresponding Jost function
$F(k)$ and briefly outline the steps to recover $V$ from $F(k).$
In Section~3 we show that the transmission eigenvalues are related to the zeros of
the key quantity $D(k)$ defined in (3.1), and in (3.4) we express $D(k)$ in terms of
the ``perturbed" Jost function $F(k)$ and the ``unperturbed" Jost function $F_0(k),$
and in (3.5) we express $D(k)$ in terms of
the ``perturbed" scattering matrix $S(k)$ and the ``unperturbed" scattering matrix $S_0(k).$ With the help of (3.5) we prove that any transmission eigenvalue $\lambda$
comes from a $k$-value related to the solution of the equation $S(k)=S_0(k)$ with $\lambda:=k^2,$ and hence we
provide a physical interpretation
of transmission eigenvalues. In Section~3 we also present various
properties of $D(k)$ in preparation for the solution of the inverse problem. In Section~4 we analyze the inverse problem of recovery
of the potential $V$ from $\cot\theta$ and
the key quantity $D(k),$ and we provide a procedure for the unique reconstruction of
$V.$
As seen from (3.7), knowledge of $D(k)$ is equivalent
to knowledge of all transmission eigenvalues
(including their multiplicities) and the constant $\gamma$ appearing in (3.8).
It is an open question whether the value of $\gamma$ and the value of $\cot\theta$
may be contained in knowledge of transmission eigenvalues.
In Section~5 we provide an independent proof of the
uniqueness for our inverse problem, namely,
we show that, assuming the existence problem is solved, there can be only one potential
corresponding to our input data set. In Section~6 we illustrate our theoretical
results with various explicit examples, such as showing that
the zero may or may not be a transmission eigenvalue and it does not have to be a simple
transmission eigenvalue, illustrating when the key quantity $D(k)$ and
the Jost function $F(k)$ may simultaneously vanish, and showing that
the number of real transmission eigenvalues may be finite or infinite.
In Section~6 we also provide an example in which we show that the constant $\gamma$
must be included in the input data set for a unique recovery of the potential,
although the potential in the example is a Dirac delta distribution and is not quite
in class $\Cal A.$ Finally, in Section~7 we indicate how some of the result
presented in the non-Dirichlet case either remain valid in the
Dirichlet case or how they are modified.

\vskip 10 pt
\noindent {\bf 2. PRELIMINARIES}
\vskip 3 pt

In this section we introduce several quantities relevant to (1.1)
with the
non-Dirichlet selfadjoint boundary condition (1.4) for some fixed value of $\theta$
in the interval $(0,\pi).$
We refer the reader to
[5,15,16] for further properties of such quantities. Recall that the potential $V$ in (1.1) is assumed
to belong to class $\Cal A$ defined in Section~1.

The Jost solution $f(k,x)$ to (1.1) is defined as the solution satisfying
$$f(k,x)=e^{ikx},\quad f'(k,x)=ik e^{ikx},\qquad x\ge b.
\tag 2.1$$
The regular solution $\varphi(k,x)$ corresponding to (1.1) and (1.4) satisfy the boundary conditions
$$\varphi(k,0)=1,\quad \varphi'(k,0)=-\cot\theta.\tag 2.2$$
The Jost function $F(k)$ for (1.1) with the boundary condition (1.4) is defined as [5,15,16]
$$F(k):=
-i[f'(k,0)+(\cot\theta)\,f(k,0)].\tag 2.3$$
Since $f(k,x)$ and $f(-k,x)$ are both solutions to (1.1)
and they are linearly independent [5,15,16] for $k\in\bold C\setminus \{0\},$ one can write $\varphi(k,x)$
as a linear combination of $f(k,x)$ and $f(-k,x)$ as
$$\varphi(k,x)=
\displaystyle\frac{1}{2k}\left[F(k)\,f(-k,x)-F(-k)\,f(k,x)\right].\tag 2.4$$

When $V(x)\equiv 0$ in (1.1) let us use the subscript $0$ to denote the quantities
corresponding to (1.1) and (1.4). From (2.1) we see that the corresponding
Jost solution $f_0(k,x)$
is given by
$$f_0(k,x)=e^{ikx},\qquad x\in{\bold R^+},\tag 2.5$$
and the corresponding regular solution $\varphi_0(k,x)$ satisfying (2.2) is given by
$$\varphi_0(k,x)=\cos kx-\displaystyle\frac{\sin kx}{k}\,\cot \theta.
\tag 2.6$$
Using (2.5) in (2.3) we obtain the corresponding Jost function $F_0(k)$ as
$$F_0(k):=
k-i\cot\theta.\tag 2.7$$

We use ${\bold C}$ for the complex plane,
${\bold C^+}$ for the open upper-half complex plane, $ {\bold C^-}$ for the open lower-half complex plane,
${\overline{\bold C^+}}$ for ${\bold C^+}\cup{\bold R},$  and
${\overline{\bold C^-}}$ for $ {\bold C^-}\cup{\bold R}.$
A bound state for the Schr\"odinger equation (1.1)
with the boundary condition (1.4) corresponds [5,15,16] to
a square-integrable solution to (1.1) satisfying (1.4).
Let us define
$$W:=\int_0^b dy\,V(y),\tag 2.8$$
where $b$ is the constant related to the support interval of $V.$

When the potential $V$ in (1.1) belongs to class $\Cal A,$
the relevant properties of the Jost
solution $f(k,x)$ and the regular solution $\varphi(k,x)$ are
summarized in the following theorem.

\noindent {\bf Theorem 2.1} {\it
Assume that the potential $V$ belongs to class $\Cal A$ and consider the corresponding
half-line
Schr\"odinger equation (1.1) with the boundary condition (1.4) for any particular
value of $\theta\in(0,\pi).$ Let $f(k,x),$
$\varphi(k,x),$ and $F(k)$ be the corresponding Jost solution, the regular solution,
and the Jost function, appearing in (2.1), (2.2), and (2.3),
respectively. Let $W$ be the real constant given in (2.8). Then:}

\item{(a)} {\it For each fixed $x\in{\bold R^+},$
the Jost solution $f(k,x)$ is entire in $k\in{\bold C}.$}

\item{(b)} {\it As $k\to\infty$ in ${\overline{\bold C^+}}$ we have}
$$f(k,0)-1+\displaystyle\frac{W}{2ik}=o\left(\displaystyle\frac{1}{k}\right),\tag 2.9$$
$$f'(k,0)-ik+\displaystyle\frac{W}{2}=o(1).\tag 2.10$$

\item{(c)} {\it As $k\to\infty$ in ${\overline{\bold C^-}}$ we have}
$$f(k,0)-1+\displaystyle\frac{W}{2ik}=e^{2ikb}\, o\left(\displaystyle\frac{1}{k}\right),\tag 2.11$$
$$f'(k,0)-ik+\displaystyle\frac{W}{2}=e^{2ikb}\,o(1) ,\tag 2.12$$
{\it where $b$ is the constant related to the support of $V.$}

\item{(d)} {\it For each fixed $x\in{\bold R^+},$
the regular solution $\varphi(k,x)$ and its $x$-derivative $\varphi'(k,x)$
are entire in $k.$}

\item{(e)} {\it The Jost function $F(k)$ is entire in $k\in\bold C.$
Its large-$|k|$ asymptotics
is given by}
$$F(k)-k-i\left(\displaystyle\frac{W}{2}-\cot\theta\right)=o(1),\qquad
k\to\infty \text{ in } {\overline{\bold C^+}},
\tag 2.13$$
$$F(k)-k-i\left(\displaystyle\frac{W}{2}-\cot\theta\right)=e^{2ikb}\,o(1),\qquad
k\to\infty \text{ in } {\overline{\bold C^-}}.
\tag 2.14$$

\item{(f)} {\it The Jost function $F(k)$ satisfies}
$$F(-k^*)=-F(k)^*,\qquad k\in{\bold C},\tag 2.15$$
{\it where the asterisk denotes complex conjugation. Thus, the zeros of
$F(k)$ occur either on the imaginary axis  in ${\bold C}$ or in pairs at points located
symmetrically with respect to the imaginary axis.}

\item{(g)} {\it The zeros of $F(k)$ in ${\bold C^+},$ if there are any,
can only occur on the positive imaginary axis; such zeros correspond to the
bound states, they are all simple, and their
number is finite. A real zero of $F(k)$
can only occur at $k=0,$ and such a zero, if it exists, must be simple.
There may be infinitely many zeros of $F(k)$ in $ {\bold C^-},$ such zeros may be nonsimple,
and they are located either on the negative imaginary axis in ${\bold C}$
or occur in pairs symmetrically located with respect to the negative imaginary axis.}

\item{(h)} {\it $F(k)$ and $F(-k)$ cannot simultaneously vanish at any
$k$-value in ${\bold C}\setminus\{0\}.$ The case $F(0)=0$ may occur, which is known
as the exceptional case, and in that case $F(k)$ has a simple zero at $k=0.$}

\noindent PROOF: The analyticity properties stated in (a),
(d), (e), and the properties listed in
(f) and (g) are already known [5,15,16]. The asymptotics in (2.9)-(2.12)
can be obtained through iteration by exploiting the
integral representations [5,7,21] for the Jost solution
and its $x$-derivative, which are respectively given by
$$f(k,x)=e^{ikx}+\displaystyle\frac{1}{k}\displaystyle\int_x^b dy\,[\sin k(y-x)]\,
V(y)\,f(k,y),\tag 2.16$$
$$f'(k,x)=ik e^{ikx}-\int_x^b dy\,[\cos k(y-x)]\,
V(y)\,f(k,y),\tag 2.17$$
where we have used the fact that the support of $V$
is confined to the interval $(0,b).$
By iterating (2.16) and (2.17) we get (2.9)-(2.12). Using
(2.9) and (2.10) in (2.3) we obtain (2.13) and (2.14). Finally, concerning (h),
the simplicity
of a possible zero of $F(k)$ at $k=0$
is already known [5,15,16], and the so-called exceptional case
 indicates that the number of bound states may change by one under a small
 perturbation of the potential. Furthermore,
if $F(k)$ and $F(-k)$ vanished at some nonzero
$k$ in ${\bold C},$ we would then get from (2.4) that
$\varphi(k,x)\equiv 0$ for that
$k$-value, contradicting (2.2). \qed

By Theorem~2.1(g) we know that the zeros of $F(k)$
in ${\bold C^+}$ correspond to the bound states.
Let us use $N$ to denote the number of bound states,
and assume that they occur at
$k=i\beta_j$ for $j=1,\dots,N.$
Associated with each bound state, there is a positive number $m_j,$
known as the bound-state norming constant,
which is defined as [5,15,16]
$$m_j:=\displaystyle\frac{1}{\sqrt{\int_0^\infty dx\, [f(i\beta_j,x)]^2 }},
\tag 2.18$$
where $f(k,x)$ is the Jost solution to (1.1) appearing in (2.1).

The scattering matrix $S(k)$ corresponding to (1.1)
with the boundary condition (1.4) is defined [5,15,16] as
$$S(k):=
-\displaystyle\frac{F(-k)}{F(k)}
,\tag 2.19$$
where $F(k)$ is the Jost function given in (2.3).
 From (2.19) it is seen that $S(k)$ is a
complex-valued scalar quantity even though it is called a matrix
in the physics literature. Note that we suppress the dependence on the parameter
$\theta$ in our notation
for various quantities such as $\varphi(k,x),$ $F(k),$ and $S(k).$
Using (2.7) in (2.19) we see that the scattering matrix
$S_0(k)$ associated with (1.1) and (1.4) when $V(x)\equiv 0$ is defined as
$$S_0(k):=-\displaystyle\frac{F_0(-k)}{F_0(k)}
,\tag 2.20$$
and it is given by
$$S_0(k)=\displaystyle\frac{k+i\cot\theta}{k-i\cot\theta}
.\tag 2.21$$

\noindent {\bf Theorem 2.2} {\it
Assume that the potential $V$ belongs to
class $\Cal A$ and consider the corresponding Schr\"odinger equation (1.1) on the
half line with the boundary condition (1.4) for any particular
value of $\theta\in(0,\pi).$ Let $F(k)$
and $S(k)$ be the corresponding Jost function and the
scattering matrix defined in (2.3) and (2.19), respectively. Then:}

\item{(a)} {\it The scattering matrix $S(k)$ is meromorphic
in ${\bold C},$ its poles in ${\bold C^+}$ can only occur on the positive imaginary axis, and
 such poles are simple and correspond to
the bound states of (1.1) with the boundary condition (1.4).
As a consequence of the compact-support property of $V,$
the value of the norming constant defined in (2.18) corresponding
to a bound state at $k=i\beta_j$ is uniquely
determined by the residue of $S(k)$ at $k=i\beta_j$ as}
$$m_j=
\displaystyle\sqrt{-i\, {\text{\rm{Res}}}\left(S(k),i\beta_j\right)}.\tag 2.22$$

\item{(b)} {\it As $k\to\pm\infty$ in ${\bold R},$ the large-$|k|$ asymptotics of the scattering matrix
$S(k)$ is given by}
$$S(k)=1-\displaystyle\frac{iW}{k}+\displaystyle\frac{2i}{k}\,\cot\theta+o\left(\displaystyle\frac{1}{k}\right)
,$$
{\it where $W$ is the constant defined in (2.8).}

\item{(c)} {\it The potential $V$ is uniquely
determined by the corresponding scattering matrix $S(k).$ Hence,
the potential $V$ is uniquely determined also by the Jost function
$F(k).$}

\noindent PROOF: The first statement in (a) follows from (2.19), Theorem~2.1(g), and
Theorem~2.1(h).
The proof of (2.22) is similar to the proof of Proposition~5.1(f) of [4]. We obtain (b) by using (2.13) in (2.19).
The proof of (c) is obtained as follows. From (2.19) we know that $S(k)$
is uniquely determined by $F(k).$
The zeros in ${\bold C^+}$ of $F(k)$ uniquely determine
all the bound states, and the corresponding norming constants $m_j$ are all
determined via (2.22).
We can then use the Marchenko method [5,15,16] to construct the potential
$V.$ To achieve this, we first form
the Marchenko kernel [5,15,16] defined as
$$\Omega(y):=\displaystyle\frac{1}{2\pi}\int_{-\infty}^\infty
dk\,[S(k)-1]\,e^{iky}+
\displaystyle\sum_{j=1}^N m_j^2\, e^{-\beta_j y}.\tag 2.23$$
We next use $\Omega(y)$ as input
in the Marchenko
integral equation
$$K(x,y)+\Omega(x+y)+\displaystyle\int_x^\infty dz\,K(x,z)\,\Omega(z+y)=0,
\qquad 0<x<y,\tag 2.24$$
and obtain $K(x,y).$
The existence and uniqueness of $K(x,y)$ as the solution to (2.24)
are assured [5,15,16] when $V$ is in class $\Cal A.$
Once $K(x,y)$ is obtained, the potential $V$ is recovered as [5,15,16]
$$V(x)=-2\,\displaystyle\frac{d K(x,x)}{dx}.\tag 2.25$$
Thus, the proof of (c) is complete. \qed

\vskip 10 pt
\noindent {\bf 3. TRANSMISSION EIGENVALUES}
\vskip 3 pt

In this section we show that the transmission eigenvalues related to (1.1) and (1.4)
correspond to the zeros of the key quantity $D(k)$ to be introduced in (3.1). We express $D(k)$ in terms of the Jost functions $F(k)$ and $F_0(k)$ given in (2.3) and (2.7), respectively.
By further expressing $D(k)$ in terms of the scattering matrices $S(k)$ and
$S_0(k)$ defined in (2.19) and (2.21), respectively, we
clarify the meaning of transmission eigenvalues and their physical
interpretation and prove that all transmission eigenvalues are
obtained from $k$-values corresponding to solutions of the
equation $S_0(k)=S(k)$ in the complex plane.

Recall that the transmission eigenvalues related to (1.1) with the boundary condition
(1.4) correspond to the $\lambda$-values for which (1.5) has nontrivial solutions
$\psi$ and $\psi_0.$ Using $\lambda:=k^2,$ we see that any solution satisfying the first
and third lines of (1.5) must be a constant multiple of the regular solution
$\varphi(k,x)$ to (1.1) appearing in (2.2). Similarly, any solution to (1.1) satisfying the second
and fourth lines of (1.5) must be a constant multiple of $\varphi_0(k,x)$ given
in (2.6). As a result, the last line of (1.5) is equivalent to saying that
the column vector $\bmatrix \varphi_0(k,b)\\
\varphi'_0(k,b)\endbmatrix$ and the column vector $\bmatrix \varphi(k,b)\\
\varphi'(k,b)\endbmatrix$ are constant multiples of each other and hence
they are linearly dependent. Therefore,
the last line of (1.5) is in turn equivalent to having
$D(k)=0,$ where the quantity $D(k)$
is defined in terms of a matrix determinant as
$$D(k):=\left| \matrix\varphi_0(k,b)& \varphi(k,b)\\
\noalign{\medskip}
\varphi'_0(k,b)& \varphi'(k,b)\endmatrix \right|.\tag 3.1$$
Thus, we have shown that any transmission eigenvalue $\lambda$ associated with
(1.1) and (1.4) corresponds to a zero of $D(k),$ where the transmission eigenvalue
$\lambda$ and the zero $k$ are related to each other as $\lambda=k^2.$

 From (1.1) and (2.2) it follows that, for each fixed $x,$ the
regular solutions $\varphi(k,x)$ and $\varphi_0(k,x)$ are even functions of $k.$
Thus, (3.1) implies that $D(k)$ is an even function of $k$ in ${\bold C}$ and hence
$D(k)$ is actually a function of $k^2.$ Note that (1.1), (2.2), (2.6), and (3.1) imply that $D(k)$ is real valued when $k\in\bold R.$
Using (2.4) and (2.6) in (3.1) we can express $D(k)$ in terms of the Jost function
$F(k)$ appearing in (2.3). With the help of (2.1) we can evaluate (2.4) at $x=b,$ and we obtain
$$D(k)=\left| \matrix\displaystyle\frac{e^{ikb}+e^{-ikb}}{2}-\cot\theta
\,\displaystyle\frac{e^{ikb}-e^{-ikb}}{2ik}& \quad\displaystyle\frac{F(k)\,e^{-ikb}-F(-k)\,e^{ikb}}{2k}\\
\noalign{\medskip}
\displaystyle\frac{ik\left[e^{ikb}-e^{-ikb}\right]}{2}-\cot\theta
\,\displaystyle\frac{e^{ikb}+e^{-ikb}}{2}&\quad\displaystyle\frac{F(k)\,e^{-ikb}+F(-k)\,e^{ikb}}{2i}
\endmatrix \right|.\tag 3.2$$
Simplifying the right-hand side of (3.2) we get
$$D(k)=\displaystyle\frac{1}{2i}\left[ F(k)+F(-k)\right]+
\displaystyle\frac{\cot\theta}{2k}\left[ F(k)-F(-k)\right].\tag 3.3$$
In order to give a physical interpretation to the transmission eigenvalues corresponding to
(1.1) and (1.4), let us incorporate (2.7) into (3.3). From (2.7) and (3.3) we get
$$D(k)=\displaystyle\frac{1}{2ik}\left[F_0(k)\,F(-k)-F_0(-k)\,F(k)\right].\tag 3.4$$
With the help of (2.19) and (2.20) we can write (3.4) in terms of the scattering matrices
$S(k)$ and $S_0(k)$ as
$$D(k)=\displaystyle\frac{F(k)\,F_0(k)}{2ik}\left[S_0(k)-S(k)\right].\tag 3.5$$

The relevant properties of $D(k)$ are given in the following theorem.

\noindent {\bf Theorem 3.1} {\it Assume that the potential $V$ belongs to class $\Cal A.$ Corresponding to the Schr\"odinger equation
(1.1) with the boundary condition (1.4) for some
$\theta$ in the interval $(0,\pi),$ let $D(k)$ be the quantity defined in (3.1),
$F(k)$ be the Jost function defined in (2.3),
and $W$ be the constant defined in (2.8).
Then:}

\item{(a)} {\it $D(k)$ is entire in $k\in{\bold C}.$}

\item{(b)} {\it $D(k)$ is an even function of $k$ in ${\bold C},$ i.e. $D(-k)=D(k)$ for
$k\in{\bold C}.$}

\item{(c)} {\it $D(-k^\ast)=D(k)^\ast$ for $k\in{\bold C},$ and $D(k)=D(k)^\ast$ for $k\in{\bold R}.$}

\item{(d)} {\it $D(k)\equiv 0$ if and only if $V(x)\equiv 0.$}

\item{(e)} {\it $D(k)$ and $F(k)$ cannot vanish at the same
$k$-value in ${\bold C}$ with the exception of $k=i\cot\theta,$
where $\cot\theta$ is the parameter appearing in (1.4).
We have $F(i\cot\theta)=0$ if and only if $D(i\cot\theta)=0.$}

\item{(f)} {\it Unless $V(x)\equiv 0,$ the quantity $D(k)$ is
unbounded in ${\bold C},$
and its large-$|k|$ asymptotics is given by}
$$D(k)-\displaystyle\frac{W}{2}=e^{2b \big|\text{Im}[k]\big|}\,o(1),\qquad k\to\infty
\text{ in } {\bold C},\tag 3.6$$
{\it where ${\text{\rm{Im}}[k]}$ denotes the imaginary part of $k$
and $b$ is the constant related to the support of the potential $V.$}

\item{(g)} {\it $D(k)$ is an entire function of $\lambda$ with order not exceeding $1/2,$
where $\lambda:=k^2.$}

\item{(h)} {\it Unless $V(x)\equiv 0,$ the quantity $D(k)$ has infinitely many zeros
in ${\bold C}.$ The Hadamard
factorization of $D(k)$ has the form}
$$D(k)=\gamma\,k^{2d}\displaystyle\prod_{j=1}^\infty \left(1-\displaystyle\frac{k^2}{k_j^2}\right),
\tag 3.7$$
{\it where $\gamma$ is a nonzero constant, $d$ is a nonnegative integer,
and the $\pm k_j$-values correspond to the nonzero zeros of
$D(k)$ in ${\bold C}.$ The value of $\gamma$ is given by}
$$\gamma=\displaystyle\frac{D^{(2d)}(0)}{(2d)!},\tag 3.8$$
{\it where $D^{(j)}(k)$ denotes the $j$-th derivative of $D(k)$ with respect to $k.$}

\item{(i)} {\it Although $D(k)$ is in general
unbounded in ${\bold C},$ it is always bounded
when $k\in{\bold R},$ and we have}
$$D(k)-\displaystyle\frac{W}{2}=o(1),\qquad k\to\pm\infty \text{ in } {\bold R}.\tag 3.9$$

\item{(j)} {\it The improper singular integral defined as
$$Q(k):=\displaystyle\frac{1}{\pi i}\int_{-\infty}^\infty dt\,\displaystyle\frac{D(t)-W/2}{t-k},
\qquad k\in{\overline{\bold C^+}},\tag 3.10$$
exists as a Cauchy principal value. That is, when $k\in{\bold C^+}$ the quantity
$Q(k)$ is well defined
 with the interpretation of the integral in (3.10) as
 $$\int_{-\infty}^\infty:=\displaystyle\lim_{R\to +\infty}\int_{-R}^R.\tag 3.11$$
When $k\in{\bold R},$ the quantity
$Q(k)$ is well defined
 with the interpretation of the integral in (3.10) as}
$$\int_{-\infty}^\infty:=\displaystyle\lim_{R\to +\infty}\,
\displaystyle\lim_{\epsilon\to 0^+}\left(
   \int_{-R}^{k-\epsilon}+ \int^{R}_{k+\epsilon}\right).\tag 3.12$$

\item{(k)} {\it The quantity $M(k)$ defined as the improper integral
$$M(k):=\displaystyle\frac{1}{\pi i}\int_{-\infty}^\infty dt\,\displaystyle\frac{D(t)-W/2}{t-k-i0^+},
\qquad k\in{\overline{\bold C^+}},\tag 3.13$$
exists as a Cauchy principal value, i.e.
with the interpretation of the integral in (3.13) as in (3.11)
in the limit $R\to +\infty.$
The presence of $i0^+$ in (3.13) indicates that the value of the
integral for real $k$-values must be evaluated as a limit from within ${\bold C^+}.$}

\item{(l)} {\it The quantities $Q(k)$ and $M(k)$ defined in (3.10) and (3.13),
respectively, are analytic in ${\bold C^+}.$ The quantity $M(k)$
 is continuous in $k\in{\overline{\bold C^+}},$ and it is related to $Q(k)$ as}
$$\cases M(k)=Q(k),\qquad k\in{\bold C^+},\\
\noalign{\medskip}
M(k)=Q(k)+D(k)-\displaystyle\frac{W}{2},\qquad k\in{\bold R}.\endcases\tag 3.14$$

\noindent PROOF: As seen from Theorem~2.1(e), the Jost function $F(k)$ is an entire
function of $k,$ and hence from (3.2) it follows that $D(k)$ is entire in $k.$ The evenness of $D(k)$ in $k$ directly follows from (3.3), and in fact it has already been stated below
(3.1). We obtain the first fact in (c) by using (2.15) in (3.3), and the
second fact in (c) follows from (b) and the first fact in (c).
Let us prove (d). If $V(x)\equiv 0,$ then we must have
$F(k)\equiv F_0(k),$ and hence (3.4) yields $D(k)\equiv 0.$ Conversely, if $D(k)\equiv 0,$ from (3.5)
we see that $S(k)\equiv S_0(k)$ because we cannot have $F(k)\equiv 0$ or $F_0(k)\equiv 0$ due to (2.13). On the other hand, by Theorem~2.2(c) we know that $S(k)$ uniquely determines $V$ and hence $S_0(k)$ can only correspond to $V(x)\equiv 0.$ Thus, the proof of (d) is complete. For the proof
of (e) we proceed as follows. If $D(k)$ and $F(k)$ vanish at a nonzero $k$-value, then
(3.4) implies that we must have $F_0(k)=0$ at that $k$-value because we know by Theorem~2.1(h) that $F(k)$ and $F(-k)$ cannot vanish at the same nonzero $k$-value. Thus,
with the help of (2.7) we see that
the only nonzero $k$-value with $D(k)=F(k)=0$ occurs at $k=i\cot\theta$ provided
$F(i\cot\theta)=0$ already. Concerning $k=0,$ since $D(k)$ and $F(k)$ are entire in
$k,$ with the help of
(2.7), from (3.4) we get
$$D(0)=-i\,F(0)+(\cot\theta)\,\dot F(0),\tag 3.15$$
where an overdot denotes the $k$-derivative. By Theorem~2.1(h), a possible zero of
$F(k)$ at $k=0$ is simple and hence $\dot F(0)\ne 0$ if $F(0)=0.$ Then, from (3.15)
we conclude that $D(0)=F(0)=0$ if and only if $\cot\theta=0,$
confirming that $D(k)$ and $F(k)$ can only vanish when $k=i\cot\theta.$
In the trivial case $V(x)\equiv 0,$ we have $D(k)\equiv 0$ and $F(k)=F_0(k),$ and hence
$F(k)$ vanishes only at $k=i\cot\theta.$
Thus, the proof of (e)
is complete.
We prove (f) by using (2.13)
and (2.14) in (3.3).
As for the proof of (g) and (h), from (a) and (b) it follows that $D(k)$ is entire in $\lambda$ with $\lambda:=k^2;$ on the other hand, (3.6) indicates that $D(k)$ is of order $1/2$ in $\lambda.$ Thus, $D(k)$ has the Hadamard factorization as stated in (3.7).
If $D(k)$ had only a finite number of zeros in ${\bold C},$ from (3.7) we see that $D(k)$ would have
to be a polynomial in $k.$ However, (3.6) would then imply that $D(k)\equiv W/2$ and hence
$D(k)$ would be bounded in ${\bold C},$ which by (f) could happen only if $V(x)\equiv 0.$ Thus, the proofs
of (g) and (h) are complete. Notice that (i) is a consequence of (f).
Let is now prove (j). For $k\in{\bold C^+}$ there is no singularity at $t=k$ because
$t\in{\bold R}.$
For $k\in{\bold R},$ since $D(t)$ is entire in $t,$ we have
$$D(t)=D(k)+(t-k)\,\dot D(k)+O((t-k)^2),
\qquad t\to k \text{ in } {\bold C},$$
and hence the singularity at $t=k$ of the integrand in (3.10) can be handled
by using the Cauchy principle value involving $\epsilon\to 0^+$ as in (3.12).
On the other hand, as stated in (b), we have $D(t)=D(-t).$ Thus, we get
$$\int_{-R}^R dt\,\displaystyle\frac{D(t)-W/2}{t-k}=2k \int_0^R dt\,\displaystyle\frac{D(t)-W/2}{t^2-k^2},\tag 3.16$$
and hence, with the help of (3.9), we see that
 the integrand on the right-hand side in (3.16) behaves as
 $o(1/t^2)$ as $t\to +\infty$ and hence it is integrable at
$t=+\infty.$ Therefore, the integral in (3.10)
is well defined as a Cauchy principal value in the sense of (3.12).
Hence, the proof of (j) is complete. The proof of (k) is similar to the
proof of (j). Let us finally prove
(l). Using
$$\displaystyle\lim_{\epsilon\to 0^+}
   \int_{k-\epsilon}^{k+\epsilon}
   dt\,\displaystyle\frac{D(t)-W/2}{t-k-i0^+}=\pi i\left(D(k)-\displaystyle\frac{W}{2}\right),\tag 3.17$$
with the help of (3.11) and (3.12), we establish (3.14). From (3.10) we obtain the
derivative of $Q(k)$ with respect to $k$ as
$$\dot Q(k)=\displaystyle\frac{1}{\pi i}\int_{-\infty}^\infty dt\,\displaystyle\frac{D(t)-W/2}{(t-k)^2},
\qquad k\in{\bold C^+},\tag 3.18$$
which is well defined for $k\in{\bold C^+}$ because the integrand does not
have a singularity when $t$ is confined
to ${\bold R}.$ Furthermore, the integrand in (3.18) is integrable at
$t=\pm\infty$ as a result of (3.9). Thus, $Q(k)$ is analytic for $k\in{\bold C^+}.$
 From the first line of (3.14) and the analyticity of $Q(k)$ in ${\bold C^+},$
 we conclude the analyticity of $M(k)$ in ${\bold C^+}.$
The continuity of $M(k)$ for $k\in{\overline{\bold C^+}}$ follows automatically because the values of $M(k)$ for $k\in{\bold R},$ by definition, are obtained as a limit
as $k$ approaches ${\bold R}$ from within ${\bold C^+}.$ We remark that the discontinuity
of $Q(k)$ when $k$ moves from ${\bold C^+}$ to ${\bold R}$ is the result of the use of
the Cauchy principal
value and is related to (3.17).
Thus, we have completed the
proof of (l). \qed

We will use (3.5) to clarify the meaning and physical interpretation of transmission eigenvalues. In the next theorem we show that any transmission eigenvalue
(i.e. any $\lambda$-value with
$\lambda:=k^2$ for which (1.5) has nontrivial solutions
$\psi$ and $\psi_0$) comes from a $k$-value satisfying the equation
$S_0(k)=S(k).$
This is somehow a surprising result because as seen from (2.20) $S_0(k)$ is not defined
at $k=i\cot\theta$ and as seen from (2.19) $S(k)$ is not defined at a nonzero $k$-value
satisfying $F(k)=0.$
Nevertheless, when $\lambda=-\cot^2\theta,$ there is also another $k$-value, namely
$k=-i\cot\theta$ corresponding to the same transmission eigenvalue.
If $\lambda=0$ is a transmission eigenvalue, even though only $k=0$ corresponds to $\lambda=0,$ we still show that the zero transmission eigenvalue $\lambda=0$
comes from $S_0(0)=S(0).$
Thus, based on the result presented in the following theorem, we conclude that
any transmission eigenvalue $\lambda$ is related
to a $k$-value at which the unperturbed scattering matrix $S_0(k)$ and the perturbed
scattering matrix $S(k)$ are equal to each other. In the language of quantum mechanics,
since $\lambda$ has the interpretation of energy, we can equivalently
state that a transmission eigenvalue occurs at an energy at which the scattering from
the ``perturbed" system agrees with the scattering from
the ``unperturbed" system.

\noindent {\bf Theorem 3.2} {\it Assume that the potential $V$ belongs to class $\Cal A,$
and consider the transmission eigenvalues related
to (1.1) with the boundary condition (1.4). Then, any transmission eigenvalue
$\lambda$ comes from a $k$-value satisfying $S_0(k)=S(k),$ where $\lambda:=k^2.$}

\noindent PROOF: We first consider nonzero
transmission eigenvalues and then the zero transmission eigenvalue. Recall that a transmission eigenvalue corresponds to a zero
of $D(k)$ defined in (3.1), and hence from (3.5) we see that a nonzero
zero of $D(k)$ can occur at a $k$-value
where $F_0(k)=0,$ $F(k)=0,$ or $S_0(k)=S(k).$ From (2.7) we see that the only zero of
$F_0(k)$ occurs when $k=i\cot\theta.$ Furthermore, from Theorem~3.1(e) we know that
a transmission eigenvalue and a zero of $F(k)$ are simultaneously possible only when $k=i\cot\theta.$
Thus, we can conclude that any transmission eigenvalue, with a possible exception of
$\lambda=-\cot^2 \theta$ must come from a $k$-value satisfying $S_0(k)=S(k).$ Now let us
consider the specific case when $\lambda=-\cot^2 \theta$ is a transmission eigenvalue.
There are two subcases to consider, namely, the subcases $\cot\theta\ne 0$ and $\cot\theta=0.$ In the former case, i.e. if $\cot\theta\ne 0,$ from (3.4) we conclude that
we must have $F(i\cot\theta)=0,$ in which case Theorem~2.1(h) implies that $F(-i\cot\theta)\ne 0.$ Thus, corresponding to the nonzero transmission eigenvalue
$\lambda=-\cot^2\theta,$ from (3.5) we see that neither $F_0(k)$ nor $F(k)$ vanish at
$k=-i\cot\theta,$ and hence we must have $S_0(k)=S(k)$ satisfied at $k=-i\cot\theta.$
In fact, in this subcase, from (2.19) and (2.21) we get $S_0(-i\cot\theta)=0$ and
$S(-i\cot\theta)=0,$ and hence $S_0(-i\cot\theta)=S(-i\cot\theta)$ indeed holds.
Now, let us consider the second subcase, i.e. when $\cot\theta=0$ and $\lambda=0$
is a transmission eigenvalue. In this case, from (2.21) we see that $S_0(0)=1$ and
 from (3.15) we see that $F(0)=0.$ From (2.19), we have
$$S(k)=\displaystyle\frac{-F(0)+k\, \dot F (0)+o(k)}{F(0)+k\,
\dot F (0)+o(k)},\qquad k\to 0 \text{ in } {\bold C},$$
which yields
$$S(0)=\displaystyle\frac{\dot F (0)}{\dot F (0)}=1,\tag 3.19$$
which again tells us that $S_0(0)=S(0)$ holds. We remark that by Theorem~2.1(g) a zero
of $F(k)$ at $k=0$ must be a simple zero and hence $\dot F(0)\ne 0$ if $F(0)=0.$
Thus, (3.19) is valid. \qed

The next result shows that if $\lambda$ is a transmission eigenvalue of
(1.1) with the boundary condition (1.4) then $\lambda^\ast$ is also a transmission eigenvalue. Thus, the transmission eigenvalues are either real or appear in complex
conjugate pairs. Recall that $\lambda$ and $k$ are related to each other as
$\lambda:=k^2.$

\noindent {\bf Proposition 3.3} {\it Assume that the potential $V$ belongs to class $\Cal A,$
and let $D(k)$ be the quantity defined in (3.1). We have the following:}

\item{(a)} {\it If $\lambda$ is a transmission eigenvalue for the corresponding Schr\"odinger equation
(1.1) with the boundary condition (1.4), then $\lambda^\ast$ is also transmission  eigenvalue.}

\item{(b)} {\it All transmission eigenvalues can be obtained from the zeros
of $D(k)$ in the closed first quadrant of ${\bold C}.$ In particular, the zeros of $D(k)$ on the positive real axis yield the positive transmission eigenvalues, the zeros of $D(k)$ on the
positive imaginary axis yield the negative transmission eigenvalues, the zeros of
$D(k)$ in the open first quadrant yield the complex transmission eigenvalues,
and a possible zero of $D(k)$ at $k=0$ corresponds to the zero transmission eigenvalue
$\lambda=0.$}

\item{(c)} {\it Unless the constant $W$ given in (2.8) is zero, there cannot be an
 infinite number of
positive transmission eigenvalues.}

\noindent PROOF: From Theorem~3.1(c) we see that if $k$ is a zero of $D(k)$ then $-k^\ast$
is also a zero of $D(k).$ The corresponding transmission eigenvalues $k^2$ and $(k^\ast)^2$ are
complex conjugates of each other, proving (a). From Theorem~3.1(b)
and Theorem~3.1(c) it follows that a complex zero $k$ of $D(k)$ in the open first quadrant
in ${\bold C}$ yields a zero in the remaining three
quadrants and that the corresponding $k^2$ is complex. Theorem~3.1(b) implies that
a zero of $D(k)$ on the positive real axis yields a zero on the negative real axis and both $k$-values correspond to the same positive transmission eigenvalue
$\lambda,$ and that a zero of $D(k)$ on the positive imaginary axis yields a zero on the negative imaginary axis and
both $k$-values correspond to the same negative transmission eigenvalue $\lambda$ via $\lambda:=k^2 .$
Thus, (b) is proved. Finally, from (3.9) we see that the number of real zeros of $D(k)$
on the positive real axis must be finite unless $W=0,$ proving (c). \qed

\vskip 10 pt
\noindent {\bf 4. THE INVERSE PROBLEM}
\vskip 3 pt

The inverse problem associated with transmission eigenvalues
related to (1.1) and
(1.4) consists of the recovery of the potential $V$ and perhaps the boundary
parameter $\cot\theta$ from an appropriate data set
containing
the corresponding transmission eigenvalues. In this paper
we consider the inverse problem of the recovery of $V$ when our data
set consists of the transmission eigenvalues (including their multiplicities),
the boundary parameter $\cot\theta,$
and the constant $\gamma$ appearing in (3.8). As in [3,4]
we define the multiplicity of a transmission eigenvalue
$\lambda$ as the multiplicity of $k^2$ as a zero of $D(k).$
In other words, we are interested
in determining $V$ when $\cot\theta$ and the quantity $D(k)$ appearing in (3.7)
are both known. We provide the unique reconstruction for this inverse
problem by using the following steps.

\item{(a)} Given $D(k),$ we use (3.9) to determine the constant
$W$ appearing in (2.8).

\item{(b)} Next, we use (3.3) and aim to determine the
corresponding Jost function $F(k)$ from knowledge of $D(k)$ and
$\cot\theta.$ By Theorem~2.2(c) we know that
$F(k)$ uniquely determines $V$ by the Marchenko procedure outlined in
the proof of Theorem~2.2. Thus, the reconstruction of $V$
will be accomplished provided
we can recover $F(k)$ from the data set consisting of
$D(k),$ $W,$ and $\cot\theta.$

\item{(c)} Motivated by (2.13),
we define $G(k)$ as
$$G(k):=F(k)-k-i\left(\displaystyle\frac{W}{2}-\cot\theta\right).\tag 4.1$$
By Theorem~2.1(e) we know that $G(k)$ is entire and satisfies
$$G(k)=o(1),\qquad
k\to\infty \text{ in } {\overline{\bold C^+}}.
\tag 4.2$$
By (2.14) we know that $G(k)$ is unbounded in
$ {\bold C^-},$ but this is irrelevant for the solution of our inverse problem because we
need $G(k)$ only for $k\in{\overline{\bold C^+}}.$
Using (4.1) in (3.3) we obtain
$$D(k)-\displaystyle\frac{W}{2}=\displaystyle\frac{1}{2i}\left[G(k)+G(-k)\right]
+\displaystyle\frac{\cot\theta}{2k}\left[G(k)-G(-k)\right],\qquad k\in{\bold R}.
\tag 4.3$$

\item{(d)} We now view (4.3) as a Riemann-Hilbert problem where $D(k)-W/2$
corresponds to the jump on ${\bold R}$ for a sectionally analytic function.
Our goal is to write the right-hand side of (4.3) as the difference of
a ``plus" function $h_+(k)$ and a ``minus" function $h_-(k),$ i.e. to write (4.3)
in the form
$$D(k)-\displaystyle\frac{W}{2}=h_+(k)-h_-(k),\qquad k\in{\bold R}.
\tag 4.4$$
By a ``plus"
function $h_+(k)$ we mean a
function which is analytic in $k\in{\bold C^+},$ continuous for
$k\in{\overline{\bold C^+}},$ and $o(1)$ as $k\to\infty$ in ${\overline{\bold C^+}}.$ By a ``minus" function
$h_-(k)$ we mean a
function which is analytic in $k\in {\bold C^-},$ continuous for
$k\in{\overline{\bold C^-}},$ and $o(1)$ as $k\to\infty$ in ${\overline{\bold C^-}}.$ We will show that (4.4)
uniquely determines $h_+(k)$ and $h_-(k)$ when the potential $V$ belongs to class $\Cal A.$

\item{(e)}
With the help of the constant
$G(0),$ which, by (4.1), is given as
$$G(0)=F(0)-i\left(\displaystyle\frac{W}{2}-\cot\theta\right),$$
we rewrite (4.3) as
$$D(k)-\displaystyle\frac{W}{2}=\displaystyle\frac{G(k)}{2i}+
\displaystyle\frac{\cot\theta}{2k}\left[G(k)-G(0)\right]+
\displaystyle\frac{G(-k)}{2i}-
\displaystyle\frac{\cot\theta}{2k}\left[G(-k)-G(0)\right],\qquad k\in{\bold R},
$$
or equivalently as
$$D(k)-\displaystyle\frac{W}{2}=H(k)-[-H(-k)],\qquad k\in{\bold R},
\tag 4.5$$
where we have defined
$$H(k):=\displaystyle\frac{G(k)}{2i}+
\displaystyle\frac{\cot\theta}{2k}\left[G(k)-G(0)\right].\tag 4.6$$
Because $G(k)$ is entire and satisfies (4.2), we conclude that
$H(k)$ is a ``plus" function and $-H(-k)$ is a ``minus" function
satisfying (4.4), i.e. (4.4) is satisfied by choosing
$$h_+(k)=H(k),\quad h_-(k)=-H(-k).\tag 4.7$$
Thus, we have shown that the Riemann-Hilbert problem posed in (4.4)
has a solution. Our next goal is to show that the solution is unique.

\item{(f)}
 From (4.4) and (4.5) we get
$$h_+(k)-H(k)=h_-(k)+H(-k),\qquad k\in{\bold R},$$
and hence any other ``plus" function would differ from $H(k)$ by an entire function
that is $o(1)$ as $k\to\infty$ in ${\bold C},$ and thus by Liouville's theorem we can conclude
that $H(k)$ and $-H(-k)$ are the only ``plus" and ``minus" functions,
respectively, satisfying (4.4). In fact, as seen from (4.4) and (4.7)
we can express $H(k)$ in terms
of $D(k)-W/2$ by using Plemelj's formula [12,20]
$$h_+(k)=
\displaystyle\frac{1}{2\pi i}\int_{-\infty}^\infty dt\,\displaystyle\frac{D(t)-W/2}{t-k-i0^+},
\qquad k\in{\overline{\bold C^+}},$$
where the integral is the Cauchy principal value in the sense of (3.11).
Thus, a comparison with (3.13) yields
$$H(k)=\displaystyle\frac{M(k)}{2},\tag 4.8$$
where $M(k)$ is the quantity defined in (3.13).

\item{(g)} Using (4.1) and (4.6) in (4.8), we  obtain
$$\displaystyle\frac{F(k)}{2ik}\left[k+i\cot\theta\right]-\displaystyle\frac{k}{2i}-\displaystyle\frac{W}{4}
-\displaystyle\frac{\cot\theta}{2k}\,F(0)=\displaystyle\frac{M(k)}{2}.\tag 4.9$$
Recall that $F(k)$ is entire and hence it cannot have a pole at $k=-i\cot\theta.$
Thus, evaluating (4.9) at $k=-i\cot\theta$ we get
$$\displaystyle\frac{\cot\theta}{2}-\displaystyle\frac{W}{4}-\displaystyle\frac{i\, F(0)}{2}=\displaystyle\frac{M(-i\cot\theta)}{2},$$
and hence the value of $F(0)$ is uniquely
determined by our data set consisting of $D(k)$ and $\cot\theta$ and we have
$$F(0)=i\left(\displaystyle\frac{W}{2}-\cot\theta\right)+i\, M(-i\cot\theta),\tag 4.10$$
where we recall that $M(k)$ is uniquely determined by $D(k).$
Using (4.10) in (4.9), we then recover $F(k)$ uniquely and explicitly from our data set consisting of $D(k)$ and $\cot\theta$ as
$$F(k)=\displaystyle\frac{ik}{k+i\cot\theta}\left[-ik+\displaystyle\frac{W}{2}+
M(k)+\displaystyle\frac{i\,\cot\theta}{k}\left(\displaystyle\frac{W}{2}-\cot\theta+ M(-i\cot\theta)\right)
\right].\tag 4.11$$

\item{(h)} Next, we use the Marchenko method [5,15,16] to reconstruct the potential
$V$ from $F(k)$ given in (4.11). Toward our goal,
we first use (4.11) in (2.19) and obtain the
corresponding scattering matrix $S(k).$ Since $S(k)$ is meromorphic with a finite
number of simple poles at $k=i\beta_j$
on the positive imaginary axis in ${\bold C},$ we first identify the
$\beta_j$ and the corresponding norming constants $m_j$ given in (2.22).
Then, we form the Marchenko kernel $\Omega(y)$ defined in (2.23). Finally, we
uniquely recover $V(x)$ via (2.25) from the unique solution $K(x,y)$ to the Marchenko equation
given in (2.24).

\vskip 10 pt
\noindent {\bf 5. AN INDEPENDENT PROOF OF THE UNIQUENESS}
\vskip 3 pt

Our reconstruction of the potential
$V$ provided in Section~4 from the data set consisting of $D(k)$ given in
(3.7) and the value of $\cot\theta$ in (1.4) also
establishes the uniqueness in the relevant
inverse problem. This is because the uniqueness is inherent in each step of
the reconstruction. Thus, we have already proved in Section~4 that, if
there exist two potentials $V$ and $\tilde V$ in class
$\Cal A,$ where both $V$ and $\tilde V$ correspond to the same data set consisting of
the transmission eigenvalues (including their multiplicities), the value of the constant
$\gamma$ appearing in (3.8), and the value of $\cot\theta$ appearing in (1.4),
then we must have $\tilde V(x)\equiv V(x).$ In this section, we provide an independent proof of the same
uniqueness using  the spectral theory for Sturm-Liouville operators so that additional and complementary
tools are introduced to analyze inverse problems associated with
transmission eigenvalues.

In our uniqueness proof, we need the following direct consequence of the Phragm\'en-Lindel\"of principle, which
can be found in Theorem~18.1.3 of [14].

\noindent{\bf Proposition 5.1} {\it Suppose that $g(k)$ is an entire function of finite order, and that order does not exceed
a positive constant $\rho$. Suppose $g(k)$ is bounded on a set of
rays $\arg[k] = \theta_j$ with $j = 1, 2, \dots, n$ for some positive integer $n$
in such a way that the angles between consecutive rays are less than $\pi / \rho$. Then $g(k)$ must be a
constant in the entire complex plane.}

Next we state and prove our uniqueness theorem.

\noindent {\bf Theorem 5.2} {\it Assume that there exists a potential $V$
in class $\Cal A$
corresponding to the data consisting of $D(k)$
defined in (3.1) and $\cot\theta$ appearing in (1.4). Then, $V$ must be the only
potential corresponding to the data.}

\noindent PROOF: Consider the following two boundary value problems:
$$\cases -\psi'' +V(x)\,\psi =\lambda \psi,\qquad 0<x<b,\\
\noalign{\medskip}
\psi'(0)+(\cot\theta)\,\psi(0)=0, \quad \psi(b)=0,\endcases\tag 5.1$$
$$\cases -\psi'' +V(x)\,\psi =\lambda \psi,\qquad 0<x<b,\\
\noalign{\medskip}
\psi'(0)+(\cot\theta)\,\psi(0)=0, \quad \psi'(b)=0.\endcases\tag 5.2$$
 From (2.2) and (5.1) it follows that
the eigenvalues of (5.1) correspond to the zeros of $\varphi(k,b),$ where
$\varphi(k,x)$ is the regular solution to (1.1) appearing in (2.2).
That is, if the zeros of $\varphi(k,b)$ occur at
$k=\pm \omega_j$ for $j\in\bold N,$ then the eigenvalues for
(5.1) are given by $\lambda=\omega_j^2$ for $j\in\bold N.$ Note that we use
$\bold N$ to denote the set of positive integers.
 From the
Sturm-Liouville theory it is already known
[13,22] that the eigenvalues for (5.1) are real and simple and their only
accumulation point is $+\infty.$ Similarly, the eigenvalues of (5.2) correspond to the zeros of $\varphi'(k,b),$ i.e. if the zeros of $\varphi'(k,b)$ occur at
$k=\pm \eta_j$ for $j\in\bold N,$ then the eigenvalues for
(5.2) are given by $\lambda=\eta_j^2$ for $j\in\bold N.$ It is also known [13,22]
that the eigenvalues for (5.2) are real and simple and their only
accumulation point is $+\infty.$ In fact, it is already known [13,22] that
we have the interlacing property
$$\eta_1^2<\omega_1^2<\eta_2^2<\omega_2^2<\eta_3^2<\cdots.\tag 5.3$$
To prove our uniqueness result, we will show that if $\{V,\varphi\}$ and $\{\tilde V,
\tilde\varphi\}$ correspond to the same data set $\{D,\cot\theta\},$ then we must
have $\tilde V(x)\equiv V(x).$ Note that we use $\tilde\varphi(k,x)$
to denote the regular solution satisfying (2.2) and also satisfying (1.1) but with
 $\tilde V$ instead of $V$ in (1.1). For the uniqueness, it is enough to prove
 that $\tilde \varphi(k,b)=\varphi(k,b)$ and $\tilde\varphi'(k,b)=\varphi'(k,b)$
 because it is already known [13,22] that the two spectral sets consisting of
 the zeros of $\varphi(k,b)$ and $\varphi'(k,b),$ respectively, uniquely determine
$V.$ Recall that, as a consequence of Liouville's theorem,
 an entire function vanishing at infinity
must be identically zero. Thus, it is enough to prove that
$P_1(k)$ and $P_2(k)$ are entire and they vanish as $k\to\infty$ in ${\bold C},$
where we have defined
$$P_1(k):=\displaystyle\frac{\tilde\varphi(k,b)-\varphi(k,b)}{\varphi_0(k,b)},
\quad P_2(k):=\displaystyle\frac{\tilde\varphi'(k,b)-\varphi'(k,b)}{\varphi_0'(k,b)},
\tag 5.4$$
with $\varphi_0(k,x)$ being the quantity given in (2.6).
Both the numerators and denominators in (5.4) are even functions of $k$
and we already know the simplicity of the $\lambda$-values
corresponding to the zeros of the denominators, where
$\lambda$ and $k$ are related to each other as $\lambda:=k^2.$ Thus, we are assured
that the order of a zero of each numerator in (5.4) is not less than the order
of the corresponding zero in the denominator.
Hence, from Theorem~2.1(d) it follows that $P_1(k)$ is entire
provided that $\tilde\varphi(k,b)-\varphi(k,b)=0$ whenever $\varphi_0(k,b)=0$
and that $P_2(k)$ is entire
provided that $\tilde\varphi'(k,b)-\varphi'(k,b)=0$ whenever $\varphi_0'(k,b)=0.$
Let us now show that these two provisions indeed hold.
Since $\tilde\varphi(k,x)$ and $\varphi(k,x)$ correspond to the same $D(k),$ from (3.1) we obtain
$$\left| \matrix\varphi_0(k,b)& \varphi(k,b)\\
\noalign{\medskip}
\varphi'_0(k,b)& \varphi'(k,b)\endmatrix \right|
=\left| \matrix\varphi_0(k,b)&\tilde \varphi(k,b)\\
\noalign{\medskip}
\varphi'_0(k,b)& \tilde\varphi'(k,b)\endmatrix \right|.\tag 5.5$$
Using (5.5), we get
$$\left| \matrix\varphi_0(k,b)& \tilde \varphi(k,b)-\varphi(k,b)\\
\noalign{\medskip}
\varphi'_0(k,b)& \tilde\varphi'(k,b)-\varphi'(k,b)\endmatrix \right|=0.\tag 5.6$$
 From (5.6) we see that at the zeros of $\varphi_0(k,b)$
 we must have $\tilde\varphi(k,b)-\varphi(k,b)=0$ because $\varphi_0'(k,b)$ cannot vanish
 at a zero of $\varphi_0(k,b).$ Similarly, (5.6) implies that at the zeros of $\varphi_0'(k,b)$
 we must have $\tilde\varphi'(k,b)-\varphi'(k,b)=0$ because $\varphi_0(k,b)$ cannot vanish
 at a zero of $\varphi_0'(k,b).$ Note that we have implicitly used
 (5.3), which implies that $\varphi_0(k,b)$ and
 $\varphi_0'(k,b)$ cannot vanish simultaneously.
 Having established that
$P_1(k)$ and $P_2(k)$ are entire, we will next show that they have the
$o(1/k)$-behavior as $k\to\infty$ in ${\bold C}.$ Let $F(k)$ and $\tilde F(k)$
be the Jost functions corresponding to $\{V,\varphi\}$ and $\{\tilde V,\tilde \varphi\},$
respectively, where the Jost function is defined as in (2.3), and let $W$ and
$\tilde W$ be the respective constants defined in (2.8) corresponding to $V$ and
$\tilde V,$ respectively. From (3.9) we see that $\tilde W=W$ because
we assume that $\{V,\varphi\}$ and $\{\tilde V,\tilde \varphi\}$ correspond to the
same $D(k).$ Thus, from (2.13) and (2.14) we obtain
$$\tilde F(k)-F(k)=o(1),\quad e^{2ikb}\left[\tilde F(-k)-F(-k)\right]=o(1),
\qquad k\to\infty \text { in } {\overline{\bold C^+}},\tag 5.7$$
$$\tilde F(-k)-F(-k)=o(1),\quad e^{-2ikb}\left[\tilde F(k)-F(k)\right]=o(1),
\qquad k\to\infty \text { in } {\overline{\bold C^-}}.\tag 5.8$$
Using $f(k,b)=e^{ikb}$ and $f'(k,b)=ik\, e^{ikb}$ implied by
(2.1), from (2.6) we get
$$\varphi_0(k,b)=e^{ikb}\left(\displaystyle\frac{1}{2}-\displaystyle\frac{\cot\theta}{2ik}\right)+
e^{-ikb}\left(\displaystyle\frac{1}{2}+\displaystyle\frac{\cot\theta}{2ik}\right),\tag 5.9$$
$$\varphi_0'(k,b)=ik\left[ e^{ikb}\left(\displaystyle\frac{1}{2}-\displaystyle\frac{\cot\theta}{2ik}\right)-
e^{-ikb}\left(\displaystyle\frac{1}{2}+\displaystyle\frac{\cot\theta}{2ik}\right)\right],\tag 5.10$$
and from (2.4) we obtain
$$\tilde\varphi(k,b)-\varphi(k,b)=
\displaystyle\frac{1}{2k}\left[ e^{-ikb}\left(\tilde F(k)-F(k)\right)-
e^{ikb}\left(\tilde F(-k)-F(-k)\right)\right],\tag 5.11
$$
$$\tilde\varphi'(k,b)-\varphi'(k,b)=\displaystyle\frac{-i}{2}\left[ e^{-ikb}\left(\tilde F(k)-F(k)\right)+
e^{ikb}\left(\tilde F(-k)-F(-k)\right)\right].\tag 5.12$$
Using (5.9) and (5.11) we can rewrite $P_1(k)$ defined in (5.4) in two
equivalent forms as
$$P_1(k)=\displaystyle\frac{1}{2k}\,\displaystyle\frac{\left(\tilde F(k)-F(k)\right)-
e^{2ikb}\left(\tilde F(-k)-F(-k)\right)}{e^{2ikb}\left(\displaystyle\frac{1}{2}-\displaystyle\frac{\cot\theta}{2ik}\right)+
\left(\displaystyle\frac{1}{2}+\displaystyle\frac{\cot\theta}{2ik}\right)},\tag 5.13$$
$$P_1(k)=\displaystyle\frac{1}{2k}\,\displaystyle\frac{e^{-2ikb}\left(\tilde F(k)-F(k)\right)-
\left(\tilde F(-k)-F(-k)\right)}{\left(\displaystyle\frac{1}{2}-\displaystyle\frac{\cot\theta}{2ik}\right)+
e^{-2ikb}\left(\displaystyle\frac{1}{2}+\displaystyle\frac{\cot\theta}{2ik}\right)}.\tag 5.14$$
Using (5.7) in (5.13) and (5.8) in (5.14) we get
$$P_1(k) = o\left(\displaystyle\frac{1}{k}\right), \qquad
k \to \infty \text{ in }
{\bold C}_\epsilon,\tag 5.15$$
for any $\epsilon > 0,$
where we have defined
$${\bold C}_\epsilon:= \{k \in {\bold C} \, : \, \text{Arg}[k] \in (-\pi+\epsilon,-\epsilon)\cup(\epsilon,\pi-\epsilon)\}
,$$
with $\text{Arg}[k]$ denoting the principal argument of $k,$ i.e. $\text{Arg}[k]\in(-\pi,\pi].$
The denominators of the right-hand sides of (5.13) and (5.14) have the leading terms proportional to
$(1+e^{2ikb})$ and $(1+e^{-2ikb}),$ respectively, and hence they vanish for arbitrarily large positive or negative
values of $k.$ Thus, it is not clear that the estimate $P_1(k) = o(1/k)$ holds as $k \to \pm\infty$ in
${\bold R},$ and hence it is unclear if $P_1(k)= o(1/k)$ as $k \to \infty$
in the entire complex plane ${\bold C}.$ In order to prove that $P_1(k)= o(1/k)$ indeed holds as $k \to \infty$
in ${\bold C},$ we will use Proposition~5.1.
Note that (5.15) implies that $P_1(k)$ is bounded on any rays other than the positive and
negative axes. Thus, if we can show that $P_1(k)$
is of finite order, then Proposition~5.1 guarantees that $P_1(k)$ is constant, and by (5.15) that constant must be zero.
Therefore, we only need to estimate the order of $P_1(k)$. Recall that we are using
$\lambda:=k^2.$ In view of (2.6) the quantity $\varphi_0(k,b)$ is entire in $\lambda$ with order $1 / 2$. Hence, by
the Hadamard factorization theorem we have
$$
\varphi_0(k,b) = c_1 k^{2d_1} \prod_{n=1}^{\infty} \left(1 - \frac{k^2}{a_n^2}\right),
\tag 5.16$$
where the $a_n$ are nonzero constants, $c_1$ is a nonzero constant,
 and $d_1=0$ or $d_1=1,$ as a result of the fact that
 the zeros of $\varphi_0(k,b),$ viewed as a function of $\lambda,$ are simple. The functions
$\tilde{\varphi}(k,b)$ and
$\varphi(k,b)$ are also entire in $\lambda$ with order $1 / 2$, and hence the order of $\tilde{\varphi}(k,b) - \varphi(k,b)$
as a function of $\lambda$ cannot exceed $1/2$.
Furthermore, as we have seen,
each zero of $\varphi_0(k,b)$ is also a zero of $\tilde{\varphi}(k,b) - \varphi(k,b)$. Thus, the
Hadamard factorization theorem implies that
$$
\tilde{\varphi}(k,b) - \varphi(k,b) =
c_2 k^{2d_1+2d_2} \prod_{n=1}^{\infty} \left(1 - \frac{k^2}{a_n^2}\right) \prod_{j=1}^{q} \left(1 - \frac{k^2}{b_j^2}\right),
\tag 5.17
$$
where the $b_j$ are nonzero constants, $c_2$ is a real constant, $d_2$ is a nonnegative integer, and
$q$ is either a nonnegative integer or $q=+\infty.$
In case $q = 0$  the value of
the second product in (5.17) is understood to be identically equal to $1.$ Note that
the possibility
$c_2 = 0$ is allowed. Using (5.16) and (5.17) in (5.4) we obtain
$$
P_1(k) =
\displaystyle\frac{c_2 k^{2d_2}}{c_1} \prod_{j=1}^{q} \left(1 - \frac{k^2}{b_j^2}\right).
\tag 5.18$$
With the help of Theorem~14.2.4 of [14], from (5.18) we conclude that $P_1(k)$, as a function of $\lambda$, has
order not exceeding $1$, or equivalently
the order of $P_1(k)$ as a function of $k$ cannot exceed $2$. Therefore, applying Proposition~5.1 with
 $g(k)=P_1(k)$ and $\rho = 2,$  we conclude
that $P_1(k) \equiv 0$.
In the same way it can be shown that $P_2(k) \equiv 0$.
Toward that goal, from (5.10) and (5.12) we see that
we can rewrite $P_2(k)$ defined in (5.4) in two
equivalent forms as
$$P_2(k)=\displaystyle\frac{1}{2k}\,\displaystyle\frac{\left(\tilde F(k)-F(k)\right)+
e^{2ikb}\left(\tilde F(-k)-F(-k)\right)}{-e^{2ikb}\left(\displaystyle\frac{1}{2}-\displaystyle\frac{\cot\theta}{2ik}\right)+
\left(\displaystyle\frac{1}{2}+\displaystyle\frac{\cot\theta}{2ik}\right)},\tag 5.19$$
$$P_2(k)=\displaystyle\frac{1}{2k}\,\displaystyle\frac{e^{-2ikb}\left(\tilde F(k)-F(k)\right)+
\left(\tilde F(-k)-F(-k)\right)}{-\left(\displaystyle\frac{1}{2}-\displaystyle\frac{\cot\theta}{2ik}\right)+
e^{-2ikb}\left(\displaystyle\frac{1}{2}+\displaystyle\frac{\cot\theta}{2ik}\right)}.\tag 5.20$$
Using (5.7) in (5.19) and
by using (5.8) in (5.20) we obtain $P_2(k)=o(1/k)$ as $k\to \infty$ in ${\bold C}_\epsilon,$
and with the help of Proposition~5.1 we conclude that $P_2(k)\equiv 0.$ Thus, the proof is complete. \qed

\vskip 10 pt
\noindent {\bf 6. EXAMPLES}
\vskip 3 pt

In this section we illustrate with various explicit examples
the direct and inverse problems for transmission
eigenvalues corresponding to (1.1) and (1.4).

\noindent {\bf Example 6.1} In this example, we show that
the zero transmission eigenvalue is not necessarily simple
by constructing an example with a zero transmission eigenvalue of
multiplicity two.
Let us choose the potential $V$ as
$$V(x)=\cases v,\qquad 0<x<b,\\
\noalign{\medskip}
0,\qquad x>b,\endcases\tag 6.1$$
where $v$ is a constant parameter.
By solving (1.1) we explicitly evaluate the Jost solution $f(k,x)$ satisfying
(2.1) as
$$f(k,x)=\cases\displaystyle\frac{1}{2}\left(1+\displaystyle\frac{k}{\omega}\right)
e^{i(k-\omega)b+i\omega x} +\displaystyle\frac{1}{2}\left(1-\displaystyle\frac{k}{\omega}\right)
e^{i(k+\omega)b-i\omega x},\qquad 0<x<b,\\
\noalign{\medskip}
e^{ikx},\qquad x>b,\endcases\tag 6.2$$
where we have defined
$$\omega:=\sqrt{k^2-v}.$$
The corresponding Jost function can be evaluated explicitly by using (6.2) in
(2.3) and we get
$$F(k)=e^{ikb}\left[\left(k-i
\cot\theta\right)\,\cos (\omega b)-\left(i\omega^2+k\,\cot\theta\right)\,\displaystyle\frac{\sin(\omega b)}{\omega}\right].\tag 6.3$$
The key quantity $D(k)$ defined in (3.1) is then evaluated by using (6.3) in (3.3)
and we obtain
$$\aligned D(k)=& \left(k^2+\cot^2\theta\right)\,\cos(\omega b)\,\displaystyle\frac{\sin(k b)}{k}
\\
& -
\left(\omega^2+\cot^2\theta\right)\,\cos(kb)\,\displaystyle\frac{\sin(\omega b)}{\omega}-
v\,\cot\theta\,\displaystyle\frac{\sin(\omega b)}{\omega}\,\displaystyle\frac{\sin(k b)}{k}.
\endaligned\tag 6.4$$
We know from Section~3 that a transmission eigenvalue $\lambda$ corresponds to a zero
of $D(k)$ with $\lambda:=k^2.$
One can find examples where the zero is a transmission eigenvalue and it is a transmission
eigenvalue with multiplicity $2.$ For example, for $b=1,$ $\cot\theta=0,$ and
$v=1,$ from (6.4) we obtain
$$D(k)=\sinh 1+\left(\displaystyle\frac{\cosh 1}{2}-\sinh 1\right)\,k^2+O(k^4),\qquad k\to 0
\text{ in } {\bold C},$$
and hence $\lambda=0$ is not a transmission eigenvalue.
For $b=1,$ $\cot\theta=0,$ and
$v=-\pi^2$ we have
$$D(k)=-\displaystyle\frac{k^2}{2}+O(k^4),\qquad k\to 0
\text{ in } {\bold C},$$
and hence $\lambda=0$ is a simple transmission eigenvalue.
On the other hand, $\lambda=0$ is a double transmission eigenvalue
if we choose $b=1,$ $\cot\theta=1.8818\overline{2},$ and
$v=5.8609\overline{2},$ where we use an overbar on a digit to indicate a round off.
The real transmission eigenvalues are possible; for example,
for $b=1,$ $\cot\theta=-2,$ and
$v=16\pi^2,$ we observe no zeros of $D(k)$ on the positive imaginary axis, but
six zeros on the positive real axis that are given by
$$k_1=2.4514\overline{6},\  k_2=5.5146\overline{1},
\  k_3=8.8583\overline{5},\  k_4=13.425\overline{3},\
 k_5=15.70\overline{8},\ k_6=26.777\overline{8},$$
and hence we get no negative transmission eigenvalues and
six positive transmission eigenvalues that are given by
$$\lambda_1=6.0096\overline{6},\  \lambda_2=30.41\overline{1},
\  \lambda_3=78.470\overline{4},\  \lambda_4=180.23\overline{8},\
 \lambda_5=246.7\overline{4},\ \lambda_6=717.04\overline{9}.$$
For $b=1$ and $\cot\theta=-2,$ by increasing the value of $v$ even further we observe that
the number of positive transmission eigenvalues increases. From (2.8) and (6.1) it follows that
$D(k)=vb+o(1)$ as $k\to+\infty,$ and hence by Proposition~3.3(c) we know that there cannot be infinitely many positive transmission eigenvalues. A graphical analysis of (6.4)
on the positive imaginary axis indicates that there cannot be infinitely many zeros of
$D(k)$ on the positive imaginary axis and hence the number of real transmission eigenvalues in these examples, unless $v=0$ in (6.1), is finite.

\noindent {\bf Example 6.2} In this example, we illustrate Theorem~3.1(e)
by analyzing $D(k)$ and $F(k)$ at $k=i\cot\theta,$ which is the only $k$-value
at which $D(k)$ and $F(k)$ can simultaneously vanish.
Let us use the potential given in (6.1).
 From (6.3) we see that $F(i\cot\theta)=0$ provided that
 $$v\,\displaystyle\frac{\sin\left(\sqrt{-v-\cot^2\theta} \,b\right)}{\sqrt{-v-\cot^2\theta}}=0,$$
which happens either when $v=0,$ yielding the trivial case $V(x)\equiv 0$ and $D(k)\equiv 0,$
or when the value of $v$ is given by
$$v=-\cot^2\theta-\displaystyle\frac{n^2\pi^2}{b^2},\qquad n\in\bold N,\tag 6.5$$
where we recall that $\bold N$ denotes the set of positive integers.
If (6.5) holds, then $F(i\cot\theta)$ and $D(i\cot\theta)$
are both zero.

\noindent {\bf Example 6.3} Here we provide an example
with infinitely many positive transmission eigenvalues. Because of (3.9), by choosing
a potential with $W=0,$ where $W$ is the constant appearing in (2.8), we know
that $D(k)$ must converge to zero as $k\to\pm \infty$ and hence yielding
a possibility for infinitely many
positive transmission eigenvalues. For this purpose, let us use the potential
$$V(x)=\cases v,\qquad 0<x<\displaystyle\frac{b}{2},\\
\noalign{\medskip}
-v,\qquad\displaystyle\frac{b}{2}< x<b,\\
\noalign{\medskip}
0,\qquad x>b,\endcases$$
where $v$ is a constant parameter.
For example, if we consider
the special case with
$b=1,$ $\cot\theta=0,$ and $v=1,$ we get
$$D(k)=q_1(k)-q_2(k)-q_3(k)-q_4(k),\tag 6.6$$
where we have defined
$$q_1(k):=k\,\left(\sin k\right)\,\cos \left(\displaystyle\frac{\sqrt{k^2-1}}{2}\right)
\,\cos \left(\displaystyle\frac{\sqrt{k^2+1}}{2}\right),\tag 6.7$$
$$q_2(k):=\sqrt{k^2+1}\,\left(\cos k\right)\,\cos \left(\displaystyle\frac{\sqrt{k^2-1}}{2}\right)
\,\sin \left(\displaystyle\frac{\sqrt{k^2+1}}{2}\right),\tag 6.8$$
$$q_3(k):=\sqrt{k^2-1}\,\left(\cos k\right)\,\sin \left(\displaystyle\frac{\sqrt{k^2-1}}{2}\right)
\,\cos \left(\displaystyle\frac{\sqrt{k^2+1}}{2}\right),\tag 6.9$$
$$q_4(k):=k\,\displaystyle\frac{\sqrt{k^2-1}}{\sqrt{k^2+1}}\,\left(\sin k\right)\,\sin \left(\displaystyle\frac{\sqrt{k^2-1}}{2}\right)
\,\sin \left(\displaystyle\frac{\sqrt{k^2+1}}{2}\right).\tag 6.10$$
By Theorem~3.1(c) we know that $D(k)$ is real valued when $k$ is real.
A graphical analysis of $D(k)$ for positive $k$-values indicates
that there are infinitely many zeros of $D(k)$
accumulating at $+\infty,$ and the graph of
$D(k)$ continually oscillates and asymptotically converges to zero. The graphical observation of
the existence of infinitely many positive transmission eigenvalues
can also be confirmed by determining the asymptotics of $D(k)$ as $k\to +\infty.$
With the help of the expansion
$$\sqrt{k^2\mp 1}=k\mp\displaystyle\frac{1}{2k}+O\left(\displaystyle\frac{1}{k^3}\right),\qquad
k\to+\infty,\tag 6.11$$
we obtain
$$\sin\left(\displaystyle\frac{\sqrt{k^2\mp 1}}{2}\right)=
\sin\left(\displaystyle\frac{k}{2}\right)\mp
\displaystyle\frac{1}{4k}\,
\cos\left(\displaystyle\frac{k}{2}\right)-
\displaystyle\frac{1}{32k^2}\,
\sin\left(\displaystyle\frac{k}{2}\right)
+O\left(\displaystyle\frac{1}{k^3}\right),\qquad
k\to+\infty,\tag 6.12$$
$$\cos\left(\displaystyle\frac{\sqrt{k^2\mp 1}}{2}\right)=
\cos\left(\displaystyle\frac{k}{2}\right)\pm
\displaystyle\frac{1}{4k}\,
\sin\left(\displaystyle\frac{k}{2}\right)-
\displaystyle\frac{1}{32k^2}\,
\cos\left(\displaystyle\frac{k}{2}\right)
+O\left(\displaystyle\frac{1}{k^3}\right),\qquad
k\to+\infty.\tag 6.13$$
Using (6.11)-(6.13) in (6.7)-(6.10), we obtain
the large-$k$ asymptotics of $D(k)$ given in (6.6) as
$$D(k)=\displaystyle\frac{2}{k}\,\left[
\cos\left(\displaystyle\frac{k}{2}\right)\right]\left[
\sin\left(\displaystyle\frac{k}{2}\right)\right]^3+
O\left(\displaystyle\frac{1}{k^2}\right),\qquad
k\to+\infty,\tag 6.14$$
which can also be written as as
$$D(k)=\displaystyle\frac{\sin k}{k}\left[
\sin\left(\displaystyle\frac{k}{2}\right)\right]^2+
O\left(\displaystyle\frac{1}{k^2}\right),\qquad
k\to+\infty.\tag 6.15$$
 From (6.14) and (6.15)
we conclude that
$D(k)$ converges to zero as $k\to+\infty$ with infinitely many oscillations,
by changing signs infinitely many times and hence it has infinitely many zeros
on the positive $k$-axis.
Let us denote the zeros of
$D(k)$ on the positive real axis with $k_{j-1}$ for $j\in\bold N.$
In this special case, there are no zeros of $D(k)$ on the positive imaginary axis or
at $k=0.$ The first few positive zeros of $D(k)$ are given by
$$k_0=0.55848\overline{8},\quad
k_1=3.263\overline{9},\quad
k_2=6.6838\overline{5},\quad
k_3=9.464\overline{7},\quad
k_4=12.894\overline{2},
$$
corresponding to the transmission eigenvalues
$$\lambda_0=0.3119\overline{1},\quad
\lambda_1=10.65\overline{2},\quad
\lambda_2=44.673\overline{8},\quad
\lambda_3=89.580\overline{5},\quad
\lambda_4=166.26\overline{1},
$$
where we use $\lambda_j:=k_j^2$ and 
also observe that $\lambda_j\to j^2\pi^2$ as $j\to+\infty.$

\noindent {\bf Example 6.4} Let the potential be given by
$$V(x)=c\,\delta(x-a),$$
where $c$ is a real parameter, $a$ is a positive number in the interval $(0,b),$
and $\delta(x)$ denotes the Dirac delta distribution. The corresponding Jost solution
appearing in (2.1) is given by
$$f(k,x)=\cases \left(1+\displaystyle\frac{ic}{2k}\right)e^{ikx}-\displaystyle\frac{ic}{2k}
\,e^{2ika-ikx},
\qquad 0\le x\le a,\\
\noalign{\medskip}
e^{ikx},\qquad x\ge a.\endcases\tag 6.16$$
Using (6.16) in (2.3) and (3.2), we obtain
$$D(k)=c\left[\cos(ka)-(\cot\theta)\,\displaystyle\frac{\sin(ka)}{k}\right]^2.$$
Note that the zeros of $D(k)$ are not affected by $c,$ and those zeros are determined
by $a$ and $\cot\theta$ alone.
In this case we get the value of $\gamma$ appearing in (3.8) as
$$\gamma=\cases c\left(a\,\cot\theta-1\right)^2,\qquad \cot\theta\ne\displaystyle\frac{1}{a},\\
\noalign{\medskip}
\displaystyle\frac{ca^4}{9},\qquad \cot\theta=\displaystyle\frac{1}{a},\endcases$$
and hence $\gamma$ is needed to determine $V$ uniquely;
the transmission eigenvalues alone cannot
uniquely determine the potential. Let us now analyze the zeros of $D(k)$ in the
closed first
quadrant in ${\bold C}.$ On the positive real axis there are infinitely many zeros of
$D(k)$ with multiplicity $2,$
and they are obtained by solving
$$k\,\cot(ka)=\cot\theta.$$
On the positive imaginary axis there are no zeros of $D(k).$
The small-$k$ asymptotics of $D(k)$ is given by
$$D(k)=D_0+D_2 k^2 +D_4 k^4+O(k^6),\qquad k\to 0\text { in } {\bold C},$$
where we have
$$D_0:=c\left(a\,\cot\theta-1\right)^2,\quad D_2:=-\displaystyle\frac{c a^2}{3}\,
\left(a\,\cot\theta-1\right)\left(a\,\cot\theta-3\right),$$
$$D_4:=\displaystyle\frac{2 c a^4}{45}\left[
\left(a\,\cot\theta-3\right)^2-\displaystyle\frac{3}{2}\right].$$
Thus, when  $\cot\theta=1/a,$ we get
$$D_0=0,\quad D_2=0,\quad D_4=\displaystyle\frac{c a^4}{9},$$
yielding the existence of the zero transmission eigenvalue with multiplicity $2.$

\vskip 10 pt
\noindent {\bf 7. THE DIRICHLET CASE}
\vskip 3 pt

In the previous sections we have obtained our results for (1.1) in the non-Dirichlet case, i.e. when the boundary condition is given by (1.4). In this section we briefly present
some of those results in the Dirichlet case, i.e. when the boundary condition is given by (1.3) instead of (1.4). For the analysis of the corresponding inverse problem in the Dirichlet case
and for further details, we refer
the reader to [4].

Let $f(k,x)$ be the Jost solution to (1.1) with the asymptotic condition given in (2.1).
In the Dirichlet case the Jost function is not given by (2.4), but it is given by
$f(k,0).$ The regular solution $\varphi(k,x)$ does not satisfy (2.2) but it satisfies
$$\varphi(k,0)=0,\quad \varphi'(k,0)=1,\tag 7.1$$
and it is expressed in terms of the Jost solution as
$$\varphi(k,x)=\displaystyle\frac{1}{2ik}\left[f(-k,0)\,f(k,x)-f(k,0)\,f(-k,x)\right],
\tag 7.2$$
instead of (2.4). The scattering matrix $S(k)$ is defined as
$$S(k):=\displaystyle\frac{f(-k,0)}{f(k,0)},\tag 7.3$$
instead of (2.19). The corresponding quantities for $V(x)\equiv 0,$ i.e.
the Jost function $f_0(k,0),$ the regular solution
$\varphi_0(k,x),$ and the scattering matrix $S_0(k),$ respectively, are given by
$$f_0(k,0)\equiv 1,\quad \varphi_0(k,x)=\displaystyle\frac{\sin kx}{k},\quad S_0(k)\equiv 1,
\tag 7.4$$
instead of (2.6), (2.7), and (2.20), respectively. The definition of
the key quantity $D(k)$ given in (3.1) holds also in the Dirichlet case,
but (3.4) and (3.5) are modified and are respectively obtained with the help
of (7.2)-(7.4) as
$$D(k)=\displaystyle\frac{1}{2ik}\left[f(k,0)-f(-k,0)\right],\tag 7.5$$
$$D(k)=\displaystyle\frac{f(k,0)}{2ik}\left[S_0(k)-S(k)\right].\tag 7.6$$

In the Dirichlet case Theorem~2.1(g) and Theorem~2.1(h) hold verbatim if we replace
$F(k)$ there by $f(k,0).$ In particular, $f(k,0)$ and $f(-k,0)$ cannot vanish
simultaneously for $k\in{\bold C}\setminus\{0\}$ because otherwise the second
initial condition in (7.1) would not hold. In the Dirichlet case, we still have
Theorem~3.1(a)-(d) valid.

In the Dirichlet case, the analog of (2.15) is given by
$$f(-k^\ast,0)=f(k,0)^\ast,\qquad k\in{\bold C}.\tag 7.7$$
Thus, from (7.5) and (7.7) we get
$$D(-k)=D(k),\quad D(-k^\ast)=D(k)^\ast,\qquad k\in{\bold C}.\tag 7.8$$
The transmission eigenvalues are still those $\lambda$-values corresponding
to the zeros of $D(k),$ where $\lambda:=k^2.$ Thus, from (7.8) we see that
all transmission eigenvalues can be obtained from the zeros of
$D(k)$ in the closed first quadrant in the complex $k$-plane.
The positive transmission eigenvalues are obtained from the zeros of
$D(k)$ on the positive real axis, the negative transmission eigenvalues are obtained from the zeros of
$D(k)$ on the positive imaginary axis, a possible zero transmission eigenvalue
corresponds to the zero of $D(k)$ at $k=0,$ and the complex
transmission eigenvalues correspond to the zeros of $D(k)$ in the open first
quadrant. From (7.8) we also conclude that if $\lambda$ is a transmission eigenvalue
in the Dirichlet case then $\lambda^\ast$ must also be a transmission eigenvalue.
Thus, the transmission eigenvalues in the Dirichlet case, as in the non-Dirichlet case,
must be either real or they must occur in complex conjugate pairs.

In the Dirichlet case
Proposition~3.3(a) and Proposition~3.3(b) hold. However, in Proposition~3.3(c)
the possibility of infinitely many positive transmission eigenvalues holds even when
$W\ne 0.$ This is because from (2.11) and (7.5) we obtain
$$D(k)=\displaystyle\frac{W}{2k^2}+o\left(\displaystyle\frac{1}{k^2}\right),\qquad
k\to\pm \infty\text{ in }{\bold R},$$
instead of (3.9).

In the Dirichlet case Theorem~3.1(e) needs to be modified as follows.
The key quantity $D(k)$ and the Jost function $f(k,0)$ cannot vanish
simultaneously at any $k$-value in the complex $k$-plane. For nonnegative
$k$-values this follows from (7.5) and the fact that $f(k,0)$ and $f(-k,0)$
cannot vanish at the same nonzero $k$-value. If $f(k,0)$ vanishes at $k=0$ we must then
have
$\dot f(0,0)\ne 0$ because a possible zero of $f(k,0)$ at $k=0$ must be simple [7,21].
Then, if $f(0,0)=0,$ from (7.5) we obtain $D(0)=-i\,\dot f(0,0),$
and hence we must have $D(0)\ne 0,$ confirming that
$D(k)$ and $f(k,0)$ cannot vanish simultaneously even at $k=0.$

The fact that $D(k)$ and the Jost function $f(k,0)$ cannot vanish
simultaneously at any $k$-value in the complex $k$-plane yields the following
important conclusion about the transmission eigenvalues. From (7.6) it follows
that any transmission eigenvalue must come from a $k$-value for which we have
$S(k)=S_0(k).$ Thus, Theorem~3.2 holds even in the Dirichlet case, and we can conclude that
a transmission eigenvalue $\lambda$ occurs at a $k$-value when the
``perturbed" scattering and the ``unperturbed" scattering coincide.

\vskip 5 pt

\noindent {\bf{REFERENCES}}

\item{[1]} T. Aktosun,
{\it Inverse scattering for vowel articulation with frequency-domain data,}
Inverse Problems {\bf 21}, 899--914 (2005).

\item{[2]} T. Aktosun, {\it
Inverse scattering to determine the shape of a vocal tract,} in: M. A.
Dritschel (ed.), {\it The extended field of operator theory,} Birkh\"auser,
Basel, 2007, pp. 1--16.

\item{[3]}
T. Aktosun, D. Gintides, and V. G. Papanicolaou, {\it  The uniqueness in the inverse problem for transmission eigenvalues for the spherically symmetric variable-speed wave equation,}
Inverse Problems {\bf 27}, 115004 (2011).

\item{[4]}
T. Aktosun and V. G. Papanicolaou, {\it Reconstruction of the wave speed
 from transmission eigenvalues for the spherically-symmetric variable-speed wave equation,}
Inverse Problems {\bf 29}, 065007 (2013).

\item{[5]} T. Aktosun and R. Weder,
{\it Inverse spectral-scattering problem with two sets of
discrete spectra for the radial Schr\"odinger equation,}
Inverse Problems {\bf 22}, 89--114 (2006).

\item{[6]} F. Cakoni, D. Colton, and D. Gintides,
{\it The interior transmission eigenvalue problem,}
SIAM J. Math. Anal. {\bf 42}, 2912--2921 (2010).

\item{[7]} K. Chadan and P. C. Sabatier, {\it Inverse problems in quantum
scattering theory,} 2nd ed., Springer, New York, 1989.

\item{[8]}
D. Colton and R. Kress,
{\it Inverse acoustic and electromagnetic scattering theory,}
2nd ed.,  Springer,
New York, 1998.

\item{[9]}
D. Colton and Y. J. Leung, {\it Complex eigenvalues and the inverse spectral problem for transmission eigenvalues,}
Inverse Problems {\bf 29}, 104008 (2013).

\item{[10]}
D. Colton and P. Monk,
{\it The inverse scattering problem for time-harmonic acoustic waves in an inhomogeneous medium,}
Quart. J. Mech. Appl. Math. {\bf 41}, 97--125~(1988).

\item{[11]}
D. Colton, L. P\"aiv\"arinta, and J. Sylvester,
{\it The interior transmission problem,}
Inverse Probl. Imaging {\bf 1}, 13--28 (2007).

\item{[12]} F. D. Gakhov, {\it Boundary value problems,}
Pergamon Press, Oxford, 1966.

\item{[13]} G. M. L. Gladwell, {\it Inverse problems in vibration,} 2nd ed.,
Kluwer, Dordrecht, 2004.

\item{[14]} E. Hille, {\it Analytic function theory}, vol. II, Chelsea Publishing Co., New York, 1977.

\item{[15]} B. M. Levitan, {\it Inverse Sturm-Liouville problems,}
VNU Science Press, Utrecht, 1987.

\item{[16]}  V.\; A.\; Marchenko, {\it Sturm-Liouville operators and
applications,} Birk\-h\"au\-ser, Basel, 1986.

\item{[17]}
J. R. McLaughlin and P. L. Polyakov,
{\it On the
uniqueness of a spherically symmetric speed of sound from transmission
eigenvalues,} J. Differential Equations {\bf 107}, 351--382 (1994).

\item{[18]}
J. R. McLaughlin, P. L. Polyakov, and P. E. Sacks,
{\it Reconstruction of a spherically symmetric speed of sound,}
SIAM J. Appl. Math. {\bf 54}, 1203--1223 (1994).

\item{[19]}
J. R. McLaughlin, P. E. Sacks, and M. Somasundaram,
{\it Inverse scattering in acoustic media using interior transmission eigenvalues,}
in: G. Chavent, G. Papanicolaou, P. Sacks, and W. Symes (eds.),
{\it Inverse problems in wave propagation,} Springer,
New York, 1997, pp. 357--374.

\item{[20]} N. I. Muskhelishvili, {\it Singular integral equations,}
Wolters-Noordhoff Publishing, Groningen, the Netherlands, 1958.

\item{[21]} R. G. Newton, {\it
Scattering theory of waves and particles,} 2nd ed., Springer, New York, 1982.

\item{[22]} J. P\"oschel and E. Trubowitz, {\it Inverse spectral
theory,} Academic Press, Orlando, 1987.

\end